\documentclass[preprint,11pt,a4paper]{elsarticle}
\usepackage{amsmath,amsfonts,amssymb,amsthm,graphics,color}
\usepackage[a4paper,left=3.5cm,right=3cm,top=3cm,bottom=3cm]{geometry}
\usepackage[english]{babel}
\usepackage{graphicx}
\usepackage{subcaption}
\usepackage{float}
\usepackage{pdfpages}
\usepackage{subcaption}
\usepackage{url}
\usepackage{xcolor}
\newtheorem{theorem}{Theorem}[section]

\newtheorem{assumption}{Assumption}[section]

\theoremstyle{definition}
\newtheorem{exmp}{Example}[section]

\theoremstyle{definition}

\usepackage[linesnumbered,ruled,vlined]{algorithm2e}
\SetKwSty{sc}

\numberwithin{equation}{section}
\usepackage{lipsum}
\usepackage{mathtools}
\usepackage{scalefnt}
\usepackage{empheq}
\usepackage{makeidx}
\usepackage{bm}

\journal{some journal}
\date{November 11, 2025}

\begin{document}
\title{A simple predictor-corrector scheme without order reduction\\ for advection-diffusion-reaction problems}

\author[uhel]{Thi Tam Dang\corref{cor1}}
\ead{tam.dang@helsinki.fi}

\author[uibk]{Lukas Einkemmer}
\ead{lukas.einkemmer@uibk.ac.at}

\author[uibk,disc]{Alexander Ostermann}
\ead{alexander.ostermann@uibk.ac.at}

\cortext[cor1]{Corresponding author}
\address[uhel]{Department of Mathematics and Statistics, University of Helsinki, Finland}
\address[uibk]{Department of Mathematics, University of Innsbruck, Austria}
\address[disc]{Digital Science Center, University of Innsbruck, Austria}

\begin{abstract}
Treating diffusion and advection/reaction separately is an effective strategy for solving semilinear advection-diffusion-reaction equations. However, such an approach is prone to suffer from order reduction, especially in the presence of inhomogeneous Dirichlet boundary conditions. In this paper, we extend an approach of Einkemmer and Ostermann [SIAM J.~Sci.~Comput.~37, A1577-A1592, 2015] to advection-diffusion-reaction problems, where the advection and reaction terms depend nonlinearly on both the solution and its gradient. Starting from a modified splitting method, we construct a predictor-corrector scheme that avoids order reduction and significantly improves accuracy. The predictor only requires the solution of a linear diffusion equation, while the corrector is simply an explicit Euler step of an advection-reaction equation. Under appropriate regularity assumptions on the exact solution, we rigorously establish second-order convergence for this scheme. Numerical experiments are presented to confirm the theoretical results.
\end{abstract}


\maketitle

\section{Introduction}
\label{sec:intro}

Solving large systems of time-dependent partial differential equations requires significant computational effort and cost, in general. One attractive approach is subdividing the problem into smaller subproblems, which can often be solved much more efficiently.  Examples of such methods include splitting schemes (see, e.g., \cite{HUNDSDORFER1995191, hundsdorfer2013numerical, mclachlan_quispel_2002}) and implicit-explicit (IMEX) schemes \cite{ascher97,tabernero24,wang15}.

In this paper we consider the following advection-diffusion-reaction equation
\begin{equation}\label{intro:prob}
    \partial_{t} u(t)= \Delta u(t)+f(u(t),\nabla u(t))
\end{equation}
equipped with suitable boundary conditions. It is well known that many numerical methods suffer from order reduction for non-trivial boundary conditions. In other words, the observed convergence order is not what one would expect based on the (formal) order of consistency. This is, in particular, a problem for schemes that subdivide \eqref{intro:prob} in order to avoid the computational cost of solving a stiff nonlinear system involving both diffusion and advection/reaction.

For reaction-diffusion equations, splitting methods have been studied extensively. In the classical splitting approach, a semilinear reaction-diffusion problem such as
\begin{equation}\label{intro:dr}
    \partial_{t} u(t)= \Delta u(t)+f(u(t))
\end{equation}
is decomposed into two subproblems: a linear diffusion equation $\partial_{t} u(t)= \Delta u(t)$, which can be solved efficiently by fast Fourier techniques, multigrid methods or some other fast Poisson solver (see, e.g.~\cite{MCKENNEY1995348}), and a nonlinear reaction equation $\partial_{t} u(t)= f(u(t))$, which is solved by a suitable numerical method. After each time step, these two solutions are combined in an appropriate way to obtain the searched-for numerical approximation. Splitting methods are particularly effective when the boundary conditions are simple, such as periodic or homogeneous Dirichlet conditions. In this situation, the well-known Lie and Strang splitting schemes usually achieve first- and second-order accuracy, respectively (see \cite{doi:10.1137/140994204}). However, for inhomogeneous Dirichlet boundary conditions, Strang splitting loses its second-order accuracy and performs similarly to the first-order Lie splitting. This reduction in accuracy occurs due to a mismatch between the values of the nonlinearity at the boundary and the prescribed boundary conditions (see \cite{10.1016/j.cam.2019.02.023, doi:10.1137/19M1257081, dang2025initialboundarycorrectedsplittingmethod, doi:10.1137/140994204, doi:10.1137/16M1056250}).

In recent years, modified splitting methods have been developed to avoid the order reduction, especially for diffusion-reaction problems with nontrivial boundary conditions, where the nonlinear reaction term depends only on \(u\). In a series of works, Einkemmer and Ostermann \cite{EINKEMMER201876, doi:10.1137/140994204, doi:10.1137/16M1056250} proposed modifying the splitting process by adding correction terms to both subflows. Bertoli and Vilmart later developed a five-part modified Strang splitting method that includes corrections based on the nonlinear flow \cite{doi:10.1137/19M1257081}. More recently, Dang, Einkemmer, and Ostermann \cite{dang2025initialboundarycorrectedsplittingmethod} introduced a new approach, which modifies both the initial data and the boundary conditions to restore the second-order accuracy of the Strang splitting. In addition, Nakano et al.~proposed in   \cite{KosukeNakano2019} a modified Strang splitting method for advection-diffusion-reaction equations where the gradient $\nabla u$ enters only linearly. However, an extension of their approach to the fully nonlinear case is not straightforward.

In this work, we start with the operator splitting idea introduced in \cite{doi:10.1137/140994204}. This approach transforms the original problem so that the subproblems are compatible with the prescribed boundary conditions. In the literature, this transformation has been computed in various ways \cite{doi:10.1137/19M1257081, EINKEMMER201876}. Here, we propose using a transformation that can be computed by simply evaluating the numerical solution. This enables us to generalize the scheme to advection-diffusion-reaction equations, where the nonlinearity also depends on the gradient of the solution. As an additional benefit, we also obtain a largely simplified second-order scheme that only requires the solution of a linear diffusion equation (the predictor) and a trivial, explicit corrector involving the nonlinear advection and reaction terms. We provide a rigorous convergence analysis for this predictor-corrector scheme showing that second-order accuracy is obtained for arbitrary (including time-dependent) Dirichlet boundary conditions. We also perform a series of numerical experiments that align well with the theoretical results.

The paper is organized as follows: Section~\ref{sec: splitting } introduces the considered model problem and the proposed predictor-corrector scheme. We also state our main result, The\-orem~\ref{ThA2.2}, which establishes second-order convergence under appropriate regularity assumptions. Section~\ref{sec: Lie} provides a detailed convergence analysis for a first-order integrator that does not require a correction step. Section~\ref{sec: strang} analyzes the full predictor-corrector scheme. Finally, Section~\ref{sec: experiments} presents some numerical experiments to evaluate the performance of the proposed method.

\section{Model problem and proposed predictor-corrector scheme}
\label{sec: splitting }

In this paper, we study numerical methods for a broad class of semilinear advection-diffusion-reaction problems. The diffusion is modeled by a second-order elliptic differential operator and the advection-reaction part by a smooth nonlinear function that depends on both the solution and its gradient. We consider this problem on a bounded domain \(\Omega \subset \mathbb{R}^d\) over the time interval \( t \in [0,T] \). Specifically, we consider the following initial boundary value problem:
\begin{equation}\label{A1.1}
\begin{aligned}
    \partial_{t} u(t)&= Du(t)+f\bigl(u(t),\nabla u(t)\bigr),  \\
    u(t)|_{\partial \Omega}&= b(t), \\
    u(0)&=u_0,
\end{aligned}
\end{equation}
where $D$ denotes the diffusion operator (e.g., $D=\Delta$). We assume that the nonlinearity $f$, the boundary data $b$, and the initial data $u_0$ are sufficiently smooth. In general, the boundary conditions may be time-dependent.

Building on our previous work \cite{EINKEMMER201876, doi:10.1137/140994204}, we introduce a modified splitting approach that eliminates the discrepancy between the values of the nonlinearity at the boundary and the prescribed boundary conditions. To achieve this, we seek a smooth function $q_n$ that satisfies the same boundary conditions as the nonlinearity at time $t_n$. We subtract this term from the nonlinearity and add it to the diffusion equation. Since the nonlinearity is smooth in both arguments, we can simply choose
\begin{equation}
    q_n = f(u_n, \nabla u_n),
\end{equation}
where $u_n$ is the numerical solution at time \(t = t_n\). With this modification at hand, we propose to split \eqref{A1.1} as follows:
\begin{equation}\label{2.2}
\begin{aligned}
    \partial_{t}v_n(t) &= Dv_n(t) + f(u_n,\nabla u_n),\\
    v_n(t)|_{\partial \Omega} &= b(t),
\end{aligned}
\end{equation}
and
\begin{align}\label{2.3}
    \partial_{t} w_n(t)= f\bigl(w_n(t),\nabla w_n(t)\bigr)- f(u_n,\nabla u_n),
\end{align}
where \(v_n\) and \(w_n\) denote the solutions of the diffusion and the reaction subproblem, respectively.

The key observation is that the solution of \eqref{2.3} with the initial data $w_n(t_n) = u_n$ is simply $w_n(t) = u_n$. Therefore, Lie and Strang splittings based on \eqref{2.2} and \eqref{2.3} can be considerably simplified. Either \eqref{2.3} does not need to be solved at all, or even the forward Euler method is accurate enough for its solution.

To achieve first-order accuracy, it is sufficient to neglect solving equation \eqref{2.3}. Therefore, we only need to solve \eqref{2.2} using an appropriate first-order scheme, as described in Algorithm \ref{alg 1}.

\medskip

\begin{algorithm}[H]
\KwIn{Numerical solution $u_n$ at time $t_n$}
\KwOut{Numerical solution $u_{n+1}$ at time $t_{n+1} = t_n + \tau$}
Solve \eqref{2.2} with the initial data \( v_{n}(t_n) = u_n \) to obtain \( u_{n+1} = v_n(t_n+\tau)\).
\caption{First-order scheme for equation \eqref{A1.1}}
\label{alg 1}
\end{algorithm}

\medskip

To obtain second-order accuracy, only a single explicit Euler approximation of \eqref{2.3} (the corrector) in addition to solving \eqref{2.2} (the predictor) is required. However, this time, the linear problem \eqref{2.2} must be solved with at least second-order accuracy. This can also be interpreted as a Strang splitting scheme, where the first half-step is easily seen to not change the numerical solution. This scheme is stated in Algorithm \ref{alg 2}.

\medskip

\begin{algorithm}[H]
\KwIn{Numerical solution $u_n$ at time $t_n$}
\KwOut{Numerical solution $u_{n+1}$ at time $t_{n+1} = t_n + \tau$}
Solve \eqref{2.2} with the initial data \(v_{n}(t_n)= u_n \) to obtain \( v_{n}(t_n + \tau) \)\;
Solve \eqref{2.3} with the initial data \( w_{n}(t_n +\frac{\tau}2) = v_{n}(t_n+\tau) \) using the forward Euler method to obtain \(u_{n+1} = w_{n}(t_n+\tau)\).
\caption{Second-order predictor-corrector scheme for equation \eqref{A1.1}}
\label{alg 2}
\end{algorithm}

\medskip

For the main convergence results in Theorems~\ref{ThA2.1} and~\ref{ThA2.2}, we make the following assumption about the data of \eqref{A1.1}.
\begin{assumption}\label{as1}
Let \(\Omega\subset \mathbb{R}^d\) be a domain with a smooth boundary, let \(D\) be a strongly elliptic linear operator with smooth coefficients and consider inhomogeneous Dirichlet boundary conditions with smooth boundary data. Further assume that the nonlinearity \(f\) is sufficiently smooth, and that the exact solution \(u(t)\) of \eqref{A1.1} belongs to \(C^2([0,T], H^{r+1}(\Omega))\) for some \(r > \frac{d}{2}\).
\end{assumption}

As usual, we denote the norm of $H^r$ by $\|\cdot\|_r$. Assumption \ref{as1} implies that $f$ is locally Lipschitz continuous in a strip along the exact solution $u(t)$ with respect to both arguments. Thus, there exists a constant $L(R)$ such that
\begin{align}\label{2.4}
    \Vert f(v,\nabla v)- f(w,\nabla w)\Vert_{r}\le L \left(\Vert v-w\Vert_{r} + \Vert \nabla v- \nabla w\Vert_{r}\right),
\end{align}
for $\max\left( \Vert v-u(t)\Vert_{r}, \Vert w- u(t)\Vert_{r}, \Vert \nabla v-\nabla u(t)\Vert_{r},\Vert \nabla w-\nabla u(t)\Vert_{r}\right) \le R$.

Theorem \ref{ThA2.1} is our main convergence result for the first-order scheme given in Algorithm~\ref{alg 1}.

\begin{theorem}\label{ThA2.1}
Let Assumption~\ref{as1} be satisfied. Then, there exists a constant $\tau_0>0$ such that, for all step sizes satisfying $0<\tau \le \tau_0$, the scheme in Algorithm \ref{alg 1} satisfies the global error estimate
\begin{align}\label{2.5}
    \Vert u_{n}-u(t_{n}) \Vert_{r+1} \le C\tau , \qquad 0\le t_n=n\tau \le T,
\end{align}
where the constant $C$ depends on $T$ but is independent of $\tau$ and $n$.
\end{theorem}

The next theorem presents a second-order convergence result for the predictor-corrector scheme in Algorithm \ref{alg 2}.

\begin{theorem} \label{ThA2.2}
Assume that Assumption~\ref{as1} holds. Then, there exists a constant $\tau_0>0$ such that, for all step sizes satisfying $0<\tau \le \tau_0$, the predictor-corrector scheme in Algorithm \ref{alg 2} satisfies the global error estimate
\begin{align}\label{2.6}
    \Vert u_{n}-u(t_{n}) \Vert_{r+1} \le C\tau^2(1+ \left|  \log \tau \right| ), \qquad 0\le t_n=n\tau \le T,
\end{align}
where the constant $C$ depends on $T$ but is independent of $\tau$ and $n$.
\end{theorem}

In Sections~\ref{sec: Lie} and~\ref{sec: strang}, we carry out the convergence analysis. In particular, we provide the proofs of Theorems~\ref{ThA2.1} and~\ref{ThA2.2}, respectively.

\section{Convergence analysis for the first-order scheme}
\label{sec: Lie}

This section analyzes a global error estimate for the first-order scheme in Algorithm~\ref{alg 1} within the framework of analytic semigroups (see, e.g., \cite{henry1981geometric, pazy2012semigroups}). We will prove that this method is first-order convergent.

\subsection{Analytical framework}

Following \cite{doi:10.1137/140994204}, we reformulate problem \eqref{A1.1} to enforce homogeneous boundary conditions. Note that this transformation is only used for the convergence analysis and does not affect the numerical schemes in Section~\ref{sec: splitting }. For this purpose, we introduce a smooth function \(\widetilde{u}(t)\) that solves the elliptic problem
\begin{equation}
\begin{aligned}
    D\widetilde{u}(t)&=0, \\
    \widetilde{u}(t)|_{\partial \Omega}&=b(t).
\end{aligned}	
\end{equation}
The transformed function $\widehat{u}(t)= u(t)-\widetilde{u}(t)$ then solves the following abstract evolution equation
\begin{equation}\label{A3.5}
\begin{aligned}
    \partial_{t} \widehat{u}(t) +A \widehat{u}(t) &= f\bigl(u(t), \nabla u(t)\bigr) - \partial_{t}\widetilde{u}(t),\\
    \widehat{u}(0)&= u_0 -\widetilde{u}(0).
\end{aligned}
\end{equation}
The linear operator $A$ is given by $Av = -Dv$ for all $v$ in the domain of $A$. Under these conditions, $A$ is a sectorial operator on $H^r(\Omega)$, and therefore, $-A$ is the generator of an analytic semigroup $e^{-tA}$ (see \cite[Chap.~2]{pazy2012semigroups}). Further, we assume without loss of generality that $A$ is invertible with a bounded inverse. This can be achieved by scaling $u$ appropriately. In particular, the operator $A$ satisfies the parabolic smoothing property
\begin{align}\label{A3.4}
	\left\| A^{\alpha} e^{-t A}\right\| \le Ct^{-\alpha} , \quad  \alpha \geq 0,
\end{align}
uniformly for $t\in \left( 0, T\right]  $.

For later use, we recall the definition of the so-called $\varphi$ functions
\begin{equation}\label{def-varphi}
\varphi_0(z) = e ^z,\qquad \varphi_{k+1}(z) =\frac{\varphi_k(z) - \varphi_k(0)}z,\quad k\ge 0,
\end{equation}
which are used to define the bounded operators $\varphi_k(-\tau A)$.

Note that the exact solution of \eqref{A1.1} at time \(t_{n+1} = t_{n} + \tau\) is given by the following variation of constants formula
\begin{equation}\label{A3.10}
\begin{aligned}
	u(t_{n+1})&= \widetilde{u}(t_{n+1})+ e^{-\tau A}\widehat{u}(t_{n}) + \int_0^{\tau} e^{-(\tau-s)A}  f \bigl( u(t_{n}+s), \nabla u(t_{n}+s) \bigr)\, ds\\
	& \qquad -  \int_0^{\tau} e^{-(\tau-s)A} \partial_{t}\widetilde{u}(t_n+s)\,ds.
\end{aligned}
\end{equation}

\subsection{Local error}

We will now prove the first-order consistency of Algorithm \ref{alg 1}. To do so, we will consider one step of the numerical solution with step size $\tau$, starting at time $t_{n}$ with the initial value $u(t_{n})$ on the exact solution. We have to show that the local error is $\mathcal O(\tau^2)$. To obtain a concise notation, we will use the following abbreviation
\begin{equation}\label{A3.3}
    g(u) = f(u,\nabla u).
\end{equation}
The solution $v_n(t_{n+1})$ of \eqref{2.2} with the initial value $v_n(t_n)= u(t_n)$ provides one step of Algorithm \ref{alg 1} applied to \eqref{A1.1}. Obviously, the function $\widehat v_n(t) = v_n(t) -\widetilde u(t)$ solves the problem
\begin{equation}\label{eq:vnhat}
    \partial_{t}\widehat v_n(t) + A\widehat v_n(t) = g \bigl(u(t_n)\bigr) - \partial_{t}\widetilde{u}(t), \quad \widehat v_n(t_n) = \widehat u(t_n).
\end{equation}
With the help of \eqref{A3.10} we thus get
\begin{equation}\label{vn-hat}
\begin{aligned}
	\widehat v_n(t_{n+1}) &= e^{-\tau A} \widehat{u}(t_{n}) + \tau \varphi_1(-\tau A) g\bigl(u(t_{n})\bigr) - \int_0^{\tau}e^{-(\tau-s) A}
    \partial_{t}\widetilde{u}(t_n+s) \,ds\\
    &= \widehat u(t_{n+1})- \int_0^{\tau}e^{-(\tau-s) A} \Bigl( g\bigl(u(t_{n}+s)\bigr)- g\bigl(u(t_{n})\bigr)\Bigr)ds\\[2mm]
    &= \widehat u(t_{n+1}) + \mathcal O(\tau^2),
\end{aligned}
\end{equation}
since
$$
g\bigl(u(t_{n}+s)\bigr)- g\bigl(u(t_{n})\bigr) = \mathcal O(s).
$$
The latter is an immediate consequence of the required regularity in Assumption~\ref{as1}. Therefore,
\begin{equation}\label{A3.14}
\begin{aligned}
    \mathcal{L}_{\tau}u(t_{n})& = v_n(t_{n+1}) = \widehat v_n(t_{n+1}) + \widetilde u(t_{n+1})\\
    &= u(t_{n+1}) + \mathcal{O}(\tau^2)
\end{aligned}
\end{equation}
and consequently, the local error $d_{n+1}= \mathcal{L}_{\tau}u(t_{n})- u(t_{n+1})$ satisfies
\begin{align}\label{A3.21}
	d_{n+1}= \mathcal{O}(\tau^2).
\end{align}

With the estimation of the local error, we are now in a position to derive the global error for Algorithm \ref{alg 1} applied to \eqref{A1.1}.

\subsection{Global error}

In this subsection, we will prove that Algorithm \ref{alg 1} is convergent of order one as stated in Theorem \ref{ThA2.1}. Let $e_{n}=u_{n}-u(t_{n})$ denote the global error. Note that the global error can be expressed as
\begin{align}\label{A3.22}
	e_{n+1}=\mathcal{L}_{\tau}u_{n}-\mathcal{L}_{\tau} u(t_{n})+d_{n+1},
\end{align}	
with $d_{n+1}= \mathcal{L}_{\tau} u(t_{n})- u(t_{n+1})$ being the local error. Repeating the steps \eqref{eq:vnhat}--\eqref{A3.14} with the initial value $v_n(t_n) = u_n$ shows that
\begin{equation*}
    \mathcal{L}_{\tau}u_n = \widetilde u(t_{n+1})+e^{-\tau A} \bigl(u_n-\widetilde{u}(t_{n})\bigr) + \tau \varphi_1(-\tau A) g(u_n) - \int_0^{\tau}e^{-(\tau-s) A} \partial_{t}\widetilde {u}(t_n+s) \,ds.
\end{equation*}
Taking the difference with \eqref{A3.14} shows at once
\begin{align}\label{A3.23}
	\mathcal{L}_{\tau}u_{n}-\mathcal{L}_{\tau} u(t_{n}) = e^{- \tau A}e_{n}+ \tau \varphi_1(-\tau A) G\bigl(u_n, u(t_n)\bigr),
\end{align}
where
$$
G\bigl(u_n, u(t_n)\bigr) = g(u_{n})-g\bigl(u(t_{n})\bigr).
$$
Thus, the global error can be express as follows
\begin{align}\label{A3.25}
	e_{n+1}=e^{- \tau A}e_{n}+\tau \varphi_1(-\tau A) G\bigl(u_n, u(t_n)\bigr)+d_{n+1}.
\end{align}
Solving \eqref{A3.25} and using \( e_0 = 0 \), we obtain
\begin{align*}\label{A3.26}
    e_{n}=  \tau \sum_{k=0}^{n-1} e^{-\left( n-k-1\right) \tau A}\varphi_1(-\tau A) G\bigl(u_k, u(t_k)\bigr) + \sum_{k=0}^{n-1} e^{-(n-k-1)\tau A} d_{k+1}.
\end{align*}
Using \eqref{A3.21} and the local Lipschitz bound \eqref{2.4}, we get
\begin{equation}\label{3.20}
	\Vert e_{n} \Vert_{r}\le  C\tau \sum_{k=1}^{n-1} \Big( \Vert e_{k} \Vert_{r}+ \Vert \nabla e_{k}\Vert_{r}\Big) + C\tau.
\end{equation}
Because of the smoothness of the data, the gradient $\nabla e^{-(n-k-1)\tau A} d_{k+1}$ is still bounded by $C\tau^2$. Thus, we have
\begin{equation}\label{3.21}
\begin{aligned}
    \Vert \nabla e_{n}\Vert_{r} &\le  \tau \left\|  \nabla \sum_{k=0}^{n-1} e^{-\left( n-k-1\right)\tau A} \varphi_1(-\tau A) G\bigl(u_{k},u(t_{k})\bigr)\right\| _{r} + \sum_{k=0}^{n-1}  \left\| \nabla e^{-(n-k-1)\tau A}  d_{k+1}\right\| _{r}\\
    & \le C\tau \left\Vert \nabla A^{-\frac{1}{2}}\right\Vert_{r}  \sum_{k=1}^{n-1} \left\| A^{\frac{1}{2}} e^{-(n-k-1)\tau A} \varphi_1(-\tau A)\right\| _{r} \Big( \|e_{k}\|_{r}+ \| \nabla e_{k}\|_{r} \Big) +  C \tau \\
	& \le C\tau \sum_{k=1}^{n-1} t_{n-k}^{-\frac{1}{2}} \Big( \Vert e_{k} \Vert_{r}+ \Vert \nabla e_{k}\Vert_{r}\Big) + C \tau.
\end{aligned}
\end{equation}
Combining \eqref{3.20} and \eqref{3.21}, we obtain
\begin{equation}
	\| e_{n}\|_{r+1} \le C \tau \sum_{k=1}^{n-1}  \| e_{k}\|_{r+1} + C\tau \sum_{k=1}^{n-1} t_{n-k}^{-\frac{1}{2}} \| e_{k}\|_{r+1} + C \tau.
\end{equation}
By applying a discrete Gronwall lemma, we obtain the desired bound \eqref{2.5}. This concludes the proof of Theorem~\ref{ThA2.1}.

\section{Convergence analysis for the second-order predictor-corrector scheme}
\label{sec: strang}

In this section, we perform the convergence analysis for the predictor-corrector scheme in Algorithm \ref{alg 2} applied to \eqref{A1.1}. In particular, we prove that this scheme has a global error of order two.

\subsection{Local error}

To study the local error, we carry out one step of Algorithm \ref{alg 2} of size $\tau$, starting at time $t_n$ on the exact solution $u(t_n)$. We recall that the predictor step is identical to the first-order scheme analyzed in the previous section (see \eqref{vn-hat}), i.e.,
\begin{equation}
\begin{aligned}
    \widehat v_n(t_{n+1}) &= e^{-\tau A} \widehat{u}(t_{n}) + \tau \varphi_1(-\tau A) g\bigl(u(t_{n})\bigr) - \int_0^{\tau}e^{-(\tau-s) A}
    \partial_{t}\widetilde{u}(t_n+s) \,ds\\
    v_n(t_{n+1}) &= \widetilde u(t_{n+1}) + \widehat v_n(t_{n+1}).
\end{aligned}
\end{equation}
To carry out the corrector step, we have to solve problem \eqref{2.3} on the interval $\left[ 0, \frac{\tau}{2} \right]$ with an Euler step and initial value $w_n(t_n+\frac{\tau}2) = v_n(t_{n+1})$. This yields
\begin{equation}\label{A4.5}
\begin{aligned}
	\mathcal{S}_{\tau}u(t_{n})& = \widetilde u(t_{n+1}) + \widehat v_n(t_{n+1}) + \frac\tau2\Bigl(g\bigl(v_n(t_{n+1})\bigr) - g\bigl(u(t_n)\bigr)\Bigr).
\end{aligned}
\end{equation}
Let us denote the local error by $d_{n+1}= \mathcal{S}_{\tau}u(t_{n})- u(t_{n+1})$. Taking the difference of the numerical solution \eqref{A4.5} and the exact solution \eqref{A3.10} yields the following representation of the local error
\begin{equation} \label{A4.6}
\begin{aligned}
    d_{n+1}&= \int_0^{\tau}\!\! e^{- (\tau-s) A} \bigl(g\bigl(u(t_n)\bigr) - g\bigl(u(t_{n}+s)\bigr)\bigr) ds + \frac\tau2\bigl(g\bigl(v_n(t_{n+1})\bigr) - g\bigl(u(t_n)\bigr)\bigr).
\end{aligned}
\end{equation}
Let
$$
\psi(s) = g\bigl(u(t_n+s)\bigr) = f\bigl(u(t_n+s), \nabla u(t_n+s)\bigr).
$$
Given the required regularity of Assumption \ref{as1}, we infer that
$$
\psi(s) = \psi(0) + s\psi'(0) + \mathcal O(s^2) =  g\bigl(u(t_n)\bigr) + s\, g'\bigl(u(t_n)\bigr) u'(t_n) + \mathcal O(s^2)
$$
and thus
\begin{equation} \label{A-tau3}
\begin{aligned}
\int_0^{\tau} e^{- (\tau-s) A} \bigl(\psi(0)-\psi(s)\bigr) \, ds &= -\tau^2 \varphi_2(-\tau A)g'\bigl(u(t_n)\bigr) u'(t_n) + \mathcal O(\tau^3)\\
&=-\frac12 \tau^2 g'\bigl(u(t_n)\bigr) u'(t_n) + A\mathcal{O}(\tau^3) + \mathcal{O}(\tau^3).
\end{aligned}
\end{equation}
Here, we used the identity
$$
\varphi_2(-\tau A) =\tfrac12 I -\tau A \varphi_3(-\tau A).
$$
On the other hand, we have
\begin{equation} \label{miracle}
\begin{aligned}
\frac\tau2\Bigl(g\bigl(v_n(t_{n+1})\bigr) - g\bigl(u(t_n)\bigr)\Bigr) &= \frac\tau2g'\bigl(u(t_n)\bigr)\bigl( v_n(t_{n+1}) - u(t_n)\bigr)+ \mathcal{O}(\tau^3)\\
&= \frac\tau2 g'\bigl(u(t_n)\bigr)\bigl( u(t_{n+1}) - u(t_n)\bigr)+ \mathcal{O}(\tau^3).
\end{aligned}
\end{equation}
Taking all together shows that
\begin{equation}\label{A4.17}
d_{n+1} =  A\mathcal{O}(\tau^3) + \mathcal{O}(\tau^3).
\end{equation}
This local error estimate now allows us to establish the second-order convergence result.

\subsection{Global error}

In this subsection, we will prove that the predictor-corrector scheme is second-order convergent. The proof is similar to that of the first-order scheme. The key differences are the use of the local error estimate \eqref{A4.17} and the parabolic smoothing property~\eqref{A3.4}.

Using the representation of $\mathcal{S}_{\tau}$, we obtain the error recursion
\begin{equation}\label{A4.18}
\begin{aligned}
	e_{n+1}&= \mathcal{S}_{\tau} u_{n}-  \mathcal{S}_{\tau} u({t_{n}})+d_{n+1}\\
	& = e^{- \tau A}e_{n}+\tau R(u_{n},u(t_{n}) )+ d_{n+1}+ \mathcal{O}(\tau^3),
\end{aligned}
\end{equation}
where $R$ satisfies the bound
\begin{align}
	\left\|R\bigl(u_{n},u(t_{n})\bigr)\right\|_r \le C\bigl( \|e_n\|_r + \|\nabla e_n\|_r\bigr).
\end{align}

We closely follow the proof of Theorem \ref{ThA2.1}, using the local error bound for the predictor-corrector scheme given in \eqref{A4.17}, and using the parabolic smoothing \eqref{A3.4} to bound the additional $A$ in $d_{n+1}$. Hence, we get
\begin{equation}\label{eq: 4.12}
	\Vert e_{n}\Vert_{r} \le   C\tau \sum_{k=1}^{n-1} \bigl( \Vert e_{k}\Vert_{r}+ \Vert \nabla e_{k}\Vert_{r}\bigr)+ C\tau^2(1+ \left|  \log \tau \right| ).
\end{equation}	
The gradient of $e_{n}$ can be estimated as
\begin{align}\label{eq: 4.13}
	\Vert \nabla e_{n}\Vert_{r}  \le
	C\tau \sum_{k=1}^{n-1} t_{n-k}^{-\frac{1}{2}} \bigl( \Vert e_{k} \Vert_{r}+ \Vert \nabla e_{k}\Vert_{r}\bigr)+ C\tau^2(1+ \left|  \log \tau \right| ).
\end{align}
Taking these two estimates together, we obtain
\begin{align}
	\| e_{n}\|_{r+1}&\le C\tau \sum_{k=1}^{n-1} \| e_{k}\|_{r+1}+
	C\tau \sum_{k=1}^{n-1} t_{n-k}^{-\frac{1}{2}} \| e_{k}\|_{r+1}+ C\tau^2(1+ \left|  \log \tau \right| ).
\end{align}
The proof is completed by applying a discrete Gronwall lemma.

\section{Numerical examples}
\label{sec: experiments}

This section presents numerical experiments to illustrate the convergence order of the proposed methods applied to the one-dimensional diffusion-reaction equation \eqref{A1.1} on the interval \(\Omega = (0,1)\). The diffusion operator is given by \(D = \partial_{xx}\). We consider different types of nonlinearities under various boundary conditions. The boundary values at the left and right boundaries are denoted by \(b_1\) and \(b_2\), respectively. The Laplacian is discretized using standard second-order central finite differences with \(500\) grid points. A reference solution is obtained using the classical fourth-order explicit Runge--Kutta method with sufficiently small time steps.

\begin{figure}[ht]
\centering
\begin{subfigure}{0.47\textwidth}
	\centering
	\includegraphics[width=\textwidth]{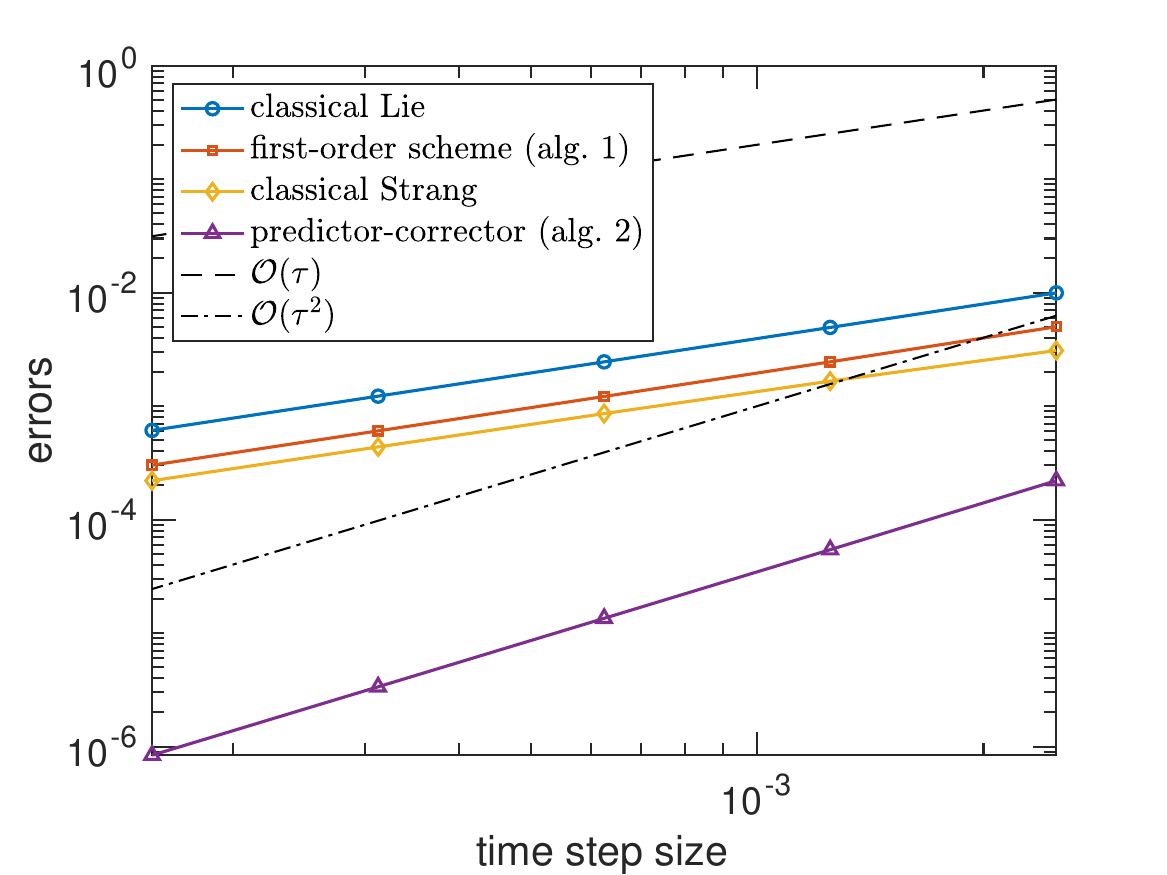}
	\caption{\(f = u (\partial_x u)^2\), time-invariant BCs}
	\label{fig1}
\end{subfigure}
\hfill
\begin{subfigure}{0.47\textwidth}
    \centering
	\includegraphics[width=\textwidth]{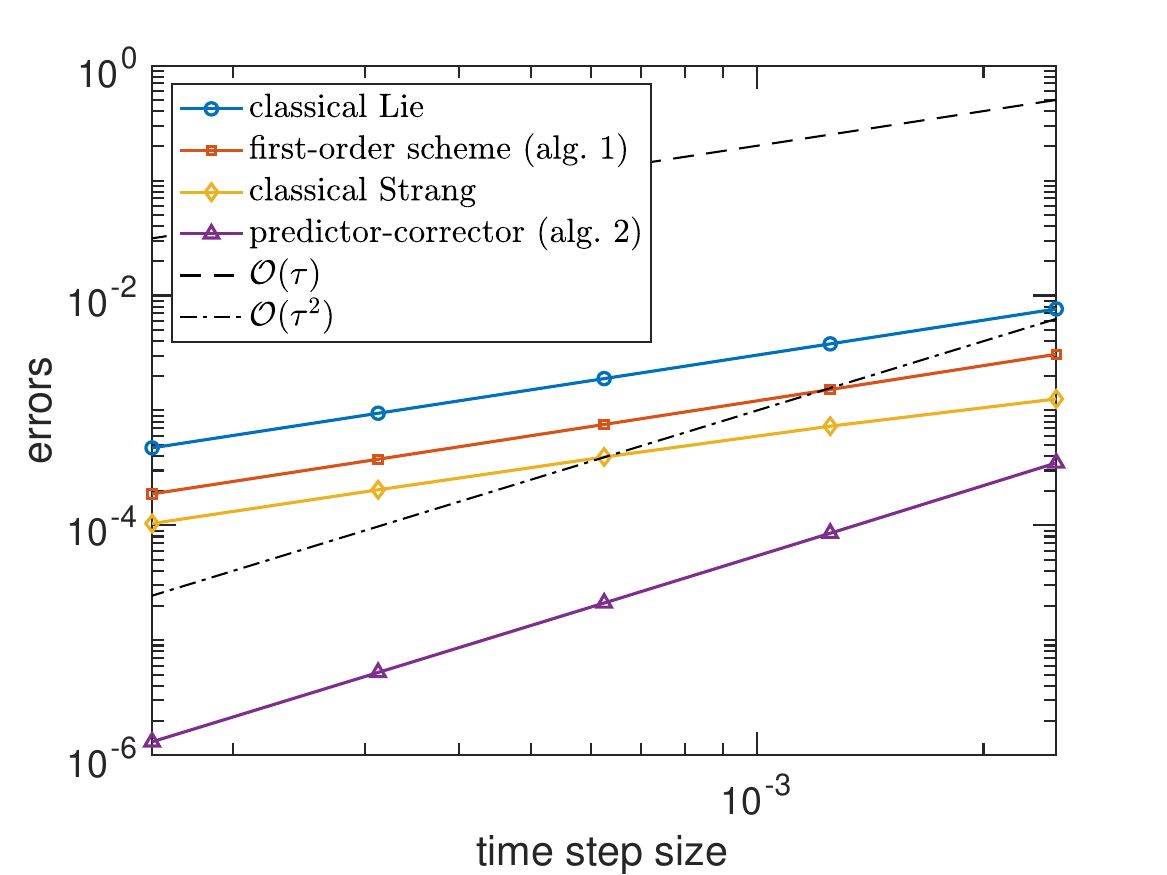}
	\caption{\(f = u^2 (\partial_x u)^2\), time-invariant BCs}
	\label{fig2}
\end{subfigure}
	
\vspace{0.5cm}
	
\begin{subfigure}{0.47\textwidth}
	\centering
	\includegraphics[width=\textwidth]{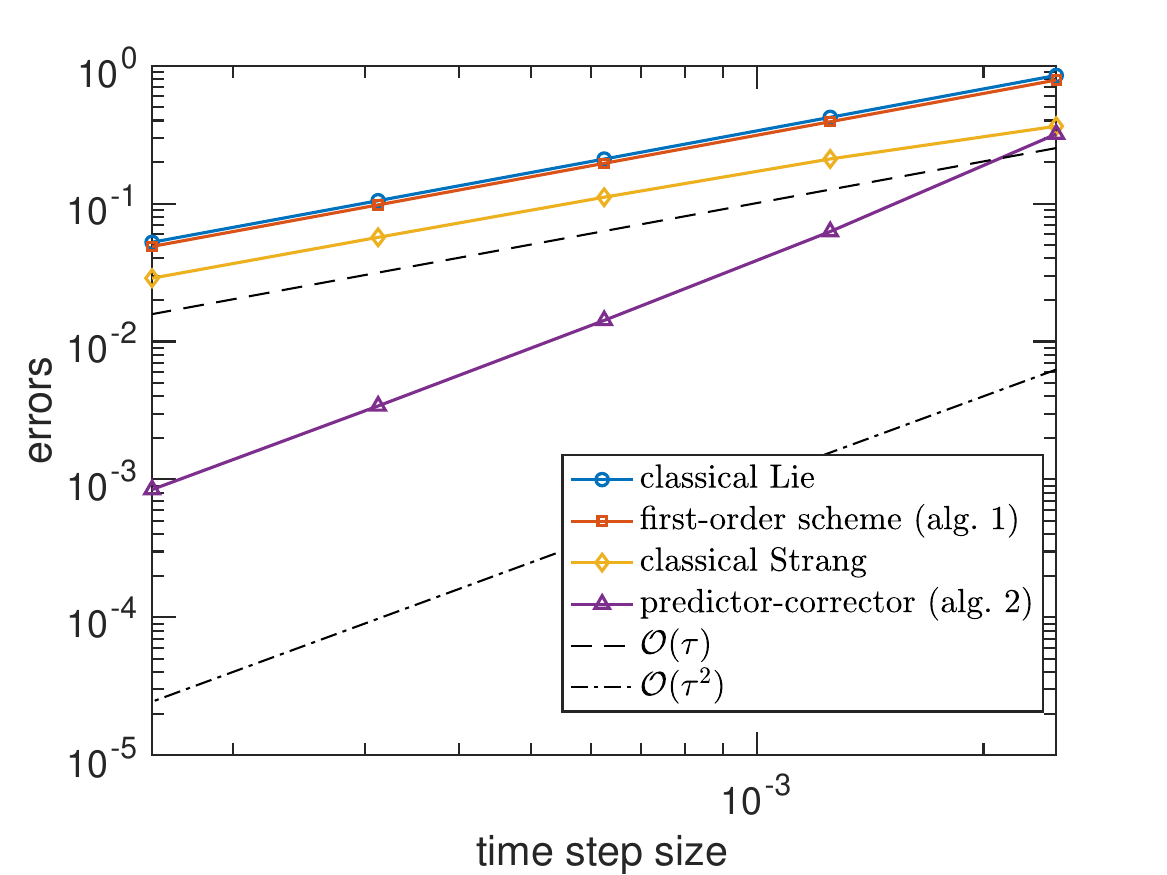}
	\caption{\(f = 3 u (\partial_x u)^2\), time-dependent BCs}
	\label{fig3}
\end{subfigure}
\hfill
\begin{subfigure}{0.47\textwidth}
	\centering
	\includegraphics[width=\textwidth]{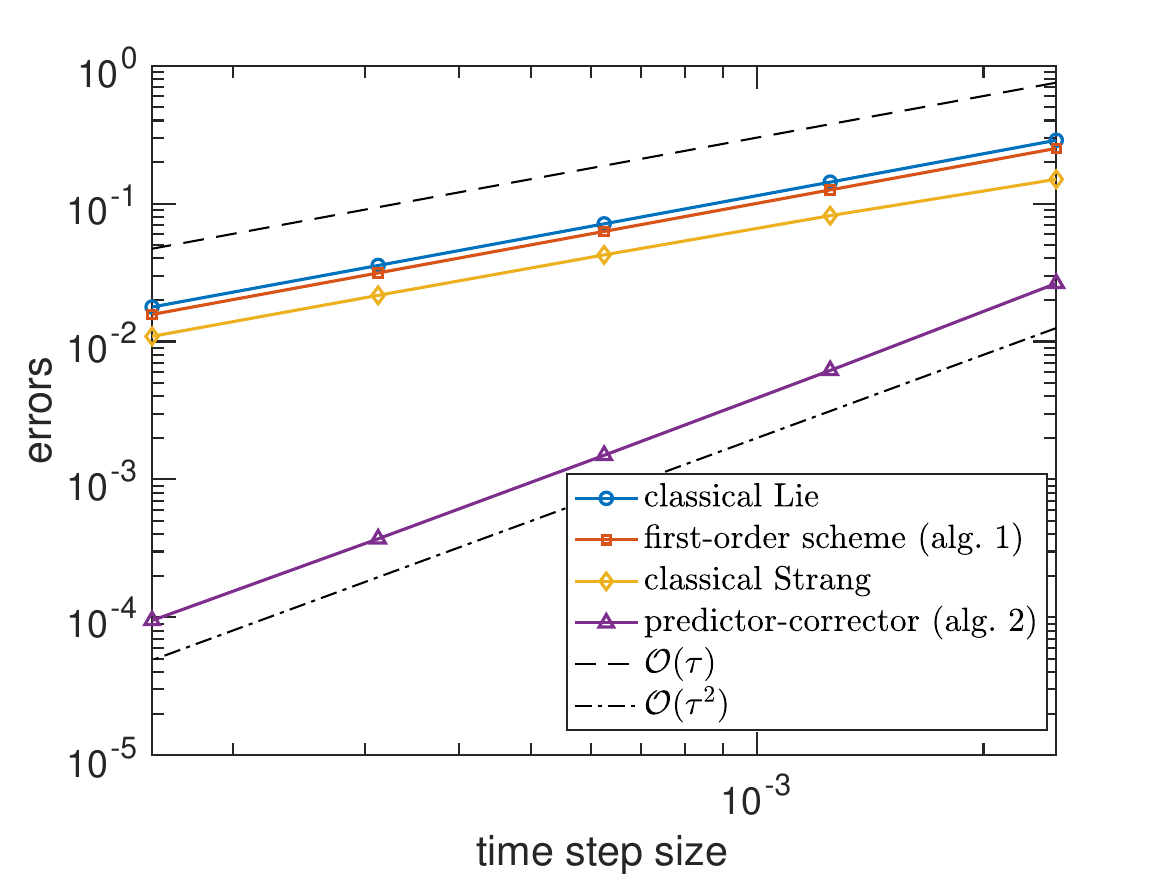}
	\caption{\(f = u(1-u) + (\partial_x u)^2\), time-dependent BCs}
	\label{fig4}
\end{subfigure}
	
\caption{We solve a one-dimensional diffusion-reaction equation with different nonlinearities and various boundary conditions. The absolute error at \( t = 0.5 \) is computed by comparing the numerical solution with a reference solution. This error is measured in a discrete \( H^{2} \) norm.}
\label{fig: main}
\end{figure}

\begin{exmp}\label{ex1}
In this example, we solve the problem \eqref{A1.1} with the nonlinear reaction term \( f(u, \partial_x u) = u(\partial_x u)^2 \). We impose time-invariant Dirichlet boundary conditions with \( b_1 = 1 \) and \( b_2 = 3 \). The initial condition is given by \( u_0(x) = 1 + 2 \sin\left( \tfrac{\pi}{2}x \right) \). The corresponding numerical results are shown in Figure~\ref{fig1}.
\end{exmp}

\begin{exmp}\label{ex2}
We consider the problem \eqref{A1.1} with the nonlinear term \( f(u, \partial_x u) = u^2(\partial_x u)^2 \), using the Dirichlet boundary conditions \( b_1 = 1 \) and \( b_2 = 1 \). The initial data is set to \(u_0 = \sin(\pi x) + 1 \), which satisfies the prescribed boundary conditions. Figure~\ref{fig2} shows the corresponding numerical results.
\end{exmp}

\begin{exmp}\label{ex3}
We study problem \eqref{A1.1} with the nonlinear reaction term \( f(u, \partial_x u) = 3u(\partial_x u)^2 \).
The left boundary condition is fixed at \( b_1 = 2 \), while the right boundary condition is time-dependent, given by \( b_2 = 1 + \cos(3 \pi t) \).
The initial condition is defined as $u_0(x) = 2 + \sin(2 \pi x)$. The numerical results are shown in Figure~\ref{fig3}.
\end{exmp}

\begin{exmp} \label{ex4}
In this example, we solve the problem \eqref{A1.1} with the nonlinear reaction term \( f(u, \partial_x u) = u(1-u)+ (\partial_x u)^2 \).  The left boundary condition is set to \( b_1 = 1 + \cos(2 \pi t)\), while the right boundary condition is given by \( b_2 = 2 + \sin(\frac\pi{2} t) \).  To satisfy these boundary conditions, we have chosen the initial data as \( u_0 = 1+ x+ \cos(2 \pi x)\). See Figure~\ref{fig4} for the numerical results.	
\end{exmp}

Figure \ref{fig: main} shows that, for the same time step size, the first-order scheme in Algorithm \ref{alg 1} is significantly more accurate than the classical Lie splitting method. This improvement is observed for both time-invariant and time-dependent boundary conditions, as well as for various nonlinearities. The figure also shows that the classical Strang splitting method reduces to approximately first-order accuracy, whereas the predictor-corrector method achieves second-order accuracy.

\section*{Acknowledgements}

This project has received funding from the European Union’s Horizon 2020 research and innovation programme under the Marie Sk\l{}odowska-Curie grant agreement No 847476.

\end{document}